\documentclass[12pt]{article}
\usepackage[utf8]{inputenc}
\usepackage{amsmath}
\usepackage{tikz-cd}
\usepackage{amssymb}
\usepackage{amsthm}
\usepackage{esvect}
\usepackage{mathtools}
\usepackage{harpoon}

\newcommand{\X}{X}

\newcommand{\defeq}{\vcentcolon=}

\newcommand{\R}{\mathbb{R}}
\newcommand\norm[1]{\left\Vert#1\right\Vert}
\newcommand\Min[1]{\mathrm{Min}(#1)}
\newcommand{\C}{\mathcal{C}}

\usepackage[shortlabels]{enumitem}
\usepackage{graphicx}
\usepackage{enumitem}
\usepackage[justification=centering]{caption}
\usepackage{subcaption}
\usepackage{comment}
\usepackage[T1]{fontenc}
\usepackage{float}
\usepackage{xcolor}

\newtheorem{theorem}{Theorem}
\newtheorem{lemma}[theorem]{Lemma}

\newtheorem*{maintheorem}{Main Theorem}
\theoremstyle{definition}
\newtheorem{example}[theorem]{Example}
\newtheorem{definition}[theorem]{Definition}
\newtheorem*{question}{Question}
\newtheorem*{acknowledge}{Acknowledgements}

\usepackage{fullpage}

\begin{document}

\begin{center}
	{\large \textsc{\textbf{A note on abelian cubulated groups}}}
\end{center}

\begin{center}
	Zachary Munro	
\end{center}

The aim of this note is to prove the following. 

\begin{maintheorem}
	An abelian group acting freely on a $\mathrm{CAT}(0)$ cube complex is free abelian.	
\end{maintheorem}

We consider actions by combinatorial automorphisms, and an action is free as long as it does not stabilize a (finite-dimensional) cube. Note that many non-free abelian groups act freely on $\mathrm{CAT}(0)$ polygonal complexes, e.g. $\mathbb{Q}$ acts freely on a $\mathrm{CAT}(0)$ 2-complex $K$ so that $\mathbb Q\backslash K$ has finite volume \cite[II.7.15 Exercises]{bridson}. 

In our following discussion, we will assume some familiarity with $\mathrm{CAT}(0)$ cube complexes. For a more detailed discussion see \cite{genevois}. $\X$ will always denote a $\mathrm{CAT}(0)$ cube complex, and $d$ the combinatorial metric on the 0-cubes $X^0$. We allow $X$ to be infinite-dimensional and even contain infinite-dimensional cubes, i.e. infinite ascending sequences of cubes. 

\begin{definition}
	A \emph{wall} $\{A, B\}$ of a set $S$ is a partition of $S$ into two nonempty subsets. An element of a wall is a \emph{halfspace}. A wall $\{A, B\}$ of a set $S$ \emph{separates} points $x,y\in S$ if $x$ and $y$ belong to different elements of the partition $\{A,B\}$. A \emph{wallspace} $Y=(S,\mathcal W)$ is a set $S$ and a set of walls $\mathcal W$ of $S$ such that any two points of $S$ are separated by at most a finite number of walls. Any subset $M\subset S$ has an \emph{induced} wallspace structure, where the walls of $M$ are the restrictions of walls of $S$. We call the dual $\mathrm{CAT}(0)$ cube complex $\C Y$ of a wallspace $Y$ the \emph{cubulation} of $Y$: The 0-cubes $\C Y^0$ correspond to choices of halfspaces for each wall such that (i) all pairs of chosen halfspaces intersect and (ii) there are no infinite descending chains $A_0\supset A_1\supset \cdots$ of chosen halfspaces. The 1-cubes $\C Y^1$ join 0-cubes with a single difference in their choice of halfspaces. The construction of $\C Y$ is completed by adding an $n$-cube wherever the 1-skeleton of an $n$-cube exists in $\C Y^1$.
	
\end{definition}

\begin{definition}
	A \emph{median algebra} $Y=(S,\mu)$ is a set $S$ with a function $\mu:S\times S\times S\to S$ such that:
	
	\begin{enumerate}
		\item $\mu(x,y,y)=y$ for all $x,y\in S$. 
		\item $\mu(x,y,z)=\mu(z,x,y)=\mu(x,z,y)$ for all $x,y,z\in S$.
		\item $\mu(\mu(x,w,y),w,z)=\mu(x,w,\mu(y,w,z))$ for all $w,x,y,z\in S$. 
	\end{enumerate}
	
	A median algebra is \emph{discrete} if for all $x,y\in S$, the set $I(x,y)\defeq\{z\in S:\mu(x,y,z)=z\}$ is finite. A set $B\subset S$ is \emph{convex} if $I(x,y)\subset B$ for all $x,y\in B$. One can give any discrete median algebra $Y$ a wallspace structure by letting the walls be all $\{B,S-B\}$ with $B$, $S-B$ convex. The \emph{cubulation} $\C Y$ of a median algebra is the cubulation of this wallspace. The 0-cells $\C Y^0$ are in bijection with the elements of $S$. A \emph{median subalgebra} is a subset $B\subset S$ closed under the $\mu$ operation. 
	
	For a $\mathrm{CAT}(0)$ cube complex $\X$, the 0-cubes $\X^0$ are naturally a discrete median algebras, where $\mu(x,y,z)$ is the unique 0-cube lying in the intersection of all halfspaces containing at least two of $x,y,z\in \X^0$. We thus think of $\mathrm{CAT}(0)$ cube complexes as discrete median algebras. We say that a hyperplane $H$ of $X$ \emph{intersects} a median subalgebra $Y\subset X$ if each halfspace of $H$ contains points of $M$.
		
	Given two median algebras $(T,\mu_T)$, $(F,\mu_F)$, the product $(T\times F, \mu_1\times\mu_2)$ defines a median algebra. We will suppress the function $\mu$ in our discussion below with the understanding that products of median algebras come equipped with the product median function. We say that a median subalgebra $Y\subset S$ is a \emph{product} of median subalgebras $T, F\subset S$ if for some (hence any) choice of basepoint $(t_1,f_1)\in T\times F$ there exists an isomorphism $T\times F\to Y\subset S$ identically mapping $(t_1,F)$ to $F$ and $(T,f_1)$ to $T$. 
		
	\end{definition}

For a $\mathrm{CAT}(0)$ cube complex $\X$, the cubulations of $X$ as a wallspace and a discrete median algebra coincide. Genevois proved that this also holds for median subalgebras of $X$ {\cite[Lemma~2.10]{genevois}}:

\begin{lemma}
\label{lem:agree}
	Let $X$ be a $\mathrm{CAT}(0)$ cube complex and $Y\subset X$ a median subalgebra. The induced wallspace structure on $Y$ coincides with its wallspace structure as a discrete median algebra.
\end{lemma}

Hyperplanes of $\mathrm{CAT}(0)$ cube complexes respect the product structure of median subalgebras. More precisely, we have the following.

\begin{lemma}
\label{lem:product}
	Let $\X$ be a $\mathrm{CAT}(0)$ cube complex. Suppose a median subalgebra $Y\subset X$ is the product $T\times F$ of median subalgebras $T,F\subset \X$. Then any hyperplane of $X$ which intersects $F$ does not intersect $T$. 
\end{lemma}

\begin{proof}
	Let $(t_1,f_1)\in Y$ be the basepoint of the product. Suppose a hyperplane $H$ intersects both $T$ and $F$. Then there must exist $(t_1,f_2)$ and $(t_2,f_1)$ separated from $(t_1,f_1)$ by $H$. 
	
	Since $M$ is a product of median algebras, $\mu((t_1,f_1),(t_1,f_2),(t_2,f_1))=(t_1,f_1)$. On the other hand, $\mu((t_1,f_1),(t_1,f_2),(t_2,f_1))$ should lie in the halfspace of $H$ containing $(t_1,f_2),(t_2,f_1)$, a contradiction.
\end{proof} 

For an isometry $g$ of a $\mathrm{CAT}(0)$ cube complex, let $\norm g=\min_{x\in X^0} d(x,gx)$ and $\Min g =\{x\in X^0 : d(x,gx)=\norm g\}$.  Haglund proved the following classification of isometries of a $\mathrm{CAT}(0)$ cube complex $\X$ \cite{haglund}.  

\begin{theorem}
\label{thm:haglund}
		For an automorphism $g$ of a $\mathrm{CAT}(0)$ cube complex $\X$, precisely one of the following holds.
		
		\begin{enumerate}
			\item $g$ is \emph{elliptic}: $g$ stabilizes a cube of $\X$. Equivalently, orbits of $g$ are bounded.
			\item $g$ is \emph{inverting}: Orbits of $g$ are unbounded and $g$ swaps the halfspaces of some hyperplane.
			\item $g$ is \emph{loxodromic}: $g$ acts as a nontrivial translation on some bi-infinite geodesic. Such a geodesic is an \emph{axis} of $g$.  
		\end{enumerate}
		
		Furthermore, when $g$ is loxodromic, every axis gets translated by $\norm g > 0$ and every element of $\Min g$ lies on an axis.
\end{theorem}

Suppose a group $G$ acts on a $\mathrm{CAT}(0)$ cube complex $\X$. Replace $\X$ by its first cubical subdivision. Then all elements $g\in G$ act by either elliptic or loxodromic isometries, and these two cases can be characterized by $\norm g=0$ and $\norm g>0$, respectively. If the action of $G$ on $\X$ is free, then all $g\in G$ act by loxodromic isometries. 

In studying isometries of cube complexes from a median viewpoint, Genevois essentially proved the following theorem \cite[Proposition~4.9]{genevois}. 
\begin{theorem}
\label{thm:ming}
		Suppose $G$ acts on a $\mathrm{CAT}(0)$ cube complex $\X$, and a central element $g\in G$ acts loxodromically. Then $\Min g$ is a $G$-invariant median subalgebra of $\X$ which decomposes as a product $T\times F$ of median subalgebras $F,T\subset X$ such that:
		
		\begin{enumerate}
			\item The action of $G$ on $\Min g=T\times F$ respects the product structure. 
			\item $F \subset X$ is an isometrically embedded $\mathrm{CAT}(0)$ cube complex that $g$ translates by $\norm g$.  
			\item $g$ acts trivially on $T$. 
		\end{enumerate}
\end{theorem}

As a result of Theorem~\ref{thm:ming}, $\Min g\subset X$ is a median subalgebra and it thus makes sense to take its cubulation $\C \Min g$. This cubulation of $\Min g$ agrees with its cubulation as an induced wallspace by Lemma~\ref{lem:agree}.

\begin{lemma}
\label{lem:commute}
	Suppose two loxodromic isometries $g,h$ 	of $\X$ commute. Then $\Min g\cap \Min {h^m}$ is nonempty for some $m>0$.
\end{lemma}

\begin{proof}
		By Theorem~\ref{thm:ming}, we have a decomposition $\Min g= T\times F$, and the induced action of $h$ on $\Min g$ respects the product structure. We consider the action of $h$ on $T$. 
		
		If the action of $h$ on $T$ is elliptic, then there is some power $h^m$ that fixes a vertex of $T$. Thus the action of $h^m$ on $\X$ stabilizes a translate of $F$ in $T\times F$. Replace $F$ with this translate. Since $F$ is an isometrically embedded $\mathrm{CAT}(0)$ cube complex and $h$ is loxodromic, there exists an axis of $h^m$ contained in $F$. 
		
		Suppose the action of $h$ on $T$ has unbounded orbits. By Theorem~\ref{thm:haglund}, after replacing $h$ by a power $h^m$ we can assume the following: $h$ acts by a loxodromic isometry on $T$, and $h$ either acts by a loxodromic isometry on $F$ or fixes a 0-cube of $F$. Let $x=(t_1,f_1)$ be an element of $\Min h$ for the action of $h$ on $\C \Min g$. Since $\Min g$ is a product, the points $t_1\in T$ and $f_1\in F$ are elements of $\Min h$ for the action of $h$ on $\C T$ and $\C F$, respectively. Thus we can find a geodesic segment $(t_1,\ldots, t_n)$ in the cubulation $\C T$ so that $\sigma_T=\bigcup_{k\in \mathbb Z} h^k\cdot (t_1,\ldots, t_n)$ is an axis of $h$ and $(t_1,\ldots, t_n)$ is a fundamental domain of $\sigma_T$. If $h$ acts by a loxodromic isometry on $F$, then there is a geodesic segment of the form $[x,hx]=(t_1,f_1),\ldots ,(t_1,f_r),(t_2,f_r),\ldots , (t_n,f_r)$, where $(f_1,\ldots, f_r)\subset F$ is a fundamental domain of an axis $\sigma_F$ of $h$ in $F$. If instead $h$ fixes a vertex of $F$, then we can take $r=1$ and $f_1$ a vertex fixed by $h$. In either case, the geodesic $\sigma = \bigcup_{k\in \mathbb Z} h^k\cdot~[x,hx]$ is an axis of $h$ for its action on $\C \Min g$. 
		
		By Theorem~\ref{thm:ming}, consecutive vertices $f_i,f_{i+1}\in F$ are distance one in $\X$, but consecutive vertices $t_i,t_{i+1}\in T$ can be distance greater than one in $\X$. For each $i=1,\ldots, n-1$, let $\gamma_i\subset \X$ be a geodesic joining $t_i$ to $t_{i+1}$. Setting $\gamma=\bigcup_{k\in \mathbb Z}h^k\cdot\big((t_1,f_1),\ldots, (t_1,f_r),(\gamma_1,f_r),\ldots, (\gamma_{n-1},f_r)\big)$, we claim that $\gamma$ is an axis of $h$ in $\X$. Indeed, suppose a hyperplane $H$ intersects $\gamma$ in two distinct edges. The two edges cannot be of the form $(h^k\cdot (t_1,f_i),h^k\cdot (t_1,f_{i+1}))$ and $(h^{k'}\cdot (t_1,f_j),h^{k'}\cdot(t_1,f_{j+1}))$ otherwise $H$ would intersect $\sigma_F$ twice. Similarly, the two edges cannot lie in geodesic segments $h^k\cdot(\gamma_i,f_r)$ and $h^{k'}\cdot(\gamma_j,f_r)$, otherwise $H$ would intersect $\sigma_T$ twice, contradicting Lemma~\ref{lem:agree}. Finally, $H$ cannot intersect an edge of the form $(h^k~\cdot~(t_1,f_i),h^k~\cdot~(t_1,f_{i+1}))$ and an edge in some $h^{k'}\cdot(\gamma_j,f_r)$ by Lemma~\ref{lem:product}. Thus no hyperplane $H$ intersects $\gamma$ twice, and $\gamma$ is a geodesic. Since $\gamma$ is stabilized by $h$, it is an axis of $h$ intersecting $F\subset \X$. 
\end{proof}

Note that in a finite-dimensional cube complex $X$, the $\mathrm{CAT}(0)$ minsets $\Min g$ and $\Min h$ intersect by the Flat Torus Theorem \cite{bridson}. However, when $X$ is infinite-dimensional, the $\mathrm{CAT}(0)$ minsets can be empty. Though the combinatorial minsets $\Min g$ and $\Min h$ are never empty, their intersection can be empty. Taking a power is required.

\begin{example}
	Consider the $\mathrm{CAT}(0)$ cube complex $\X=C\times \mathbb Z$, where $C$ is a 2-dimensional cube. Let $g$ be an isometry of $\X$ which reflects along a diagonal of $C$ and translates $\mathbb Z$ a nonzero amount. Let $h$ be an isometry of $\X$ which reflects across the other diagonal of $C$ and translates $\mathbb Z$ a nonzero amount. Then $g$, $h$ are loxodromic and $\Min g \cap \Min h$ is empty. However, $\Min g\cap \Min {h^2}$ is nonempty.
\end{example}

We now turn our attention to abelian groups. 

\begin{definition}
	Let $A$ be an abelian group. A \emph{discrete norm} on $A$ is a function $\nu:A\to \R$ such that 
	\begin{enumerate}
		\item There exists $\epsilon > 0 $ such that $\nu(g) > \epsilon$ for $g\neq 0$.
		\item $\nu(g+h)\leq \nu(g)+\nu(h)$.
		\item $\nu(g^m)=|m|\nu(g)$ for $m\in \mathbb Z$
	\end{enumerate}
\end{definition}

Stepr\=ans proved the following characterization of free abelian groups in the category of abelian groups \cite{freeAbelian}. 

\begin{theorem}
	An abelian group is free abelian if and only if it admits a discrete norm. 	
\end{theorem}


Using Stepr\=ans characterization, we prove the main theorem as follows. 

\begin{proof}[Proof of Main Theorem.]
	Let $A$ be an abelian group acting freely on a $\mathrm{CAT}(0)$ cube complex $\X$. Replace $\X$ with its cubical subdivision and define $\nu(g)\defeq \norm g$. We claim that $\nu$ is a discrete norm on $A$. Any $g\neq 0$ acts by a loxodromic isometry, so we immediately have $\nu(g)\geq 1$ and $\nu(g^m)=|m|\nu(g)$. It remains to show $\nu(g+h)\leq \nu(g)+\nu(h)$ for $g,h\in A$. 
	
	By Lemma~\ref{lem:commute}, there exists some $m>0$ such that $\Min g\cap \Min {h^m}$ is nonempty. Thus $\Min {g^m}\cap \Min {h^m}$ is also nonempty. Let $x\in \Min {g^m}\cap \Min {h^m}$. The intersection $\Min {g^m}\cap \Min {h^m}$ is stabilized by both $g^m$ and $h^m$. Hence the triangle inequality applied to $x, g^mx, h^mg^mx$ yields the following: 
	
	$$|m|\nu(gh)=\nu(h^mg^m)\leq \nu(h^m)+\nu(g^m)=|m|(\nu(g)+\nu(h)).$$
	
	Dividing by $|m|$, we retrieve our desired inequality. 
\end{proof}

The main theorem naturally leads to several questions. Namely, which solvable groups act freely on $\mathrm{CAT}(0)$ cube complexes? A construction in \cite{diadem} shows that there are arbitrarily long iterated wreath products of groups which act freely on $\mathrm{CAT}(0)$ cube complexes. Thus there are solvable groups with arbitrarily long derived length which act freely on $\mathrm{CAT}(0)$ cube complexes. In contrast, Genevois proves in \cite{genevois} that polycyclic groups acting freely on $\mathrm{CAT}(0)$ cube complexes are virtually abelian. It is thus natural to ask the following.

\begin{question}
		Which (infinitely generated) nilpotent groups act freely on $\mathrm{CAT}(0)$ cube complexes? Are they virtually free abelian?
\end{question}

\begin{acknowledge} Thank you Dani Wise for asking me the question which led to this note. You have been very encouraging along the way. Thank you Piotr Przytycki for reading through multiple drafts of this note and giving comments. This note is far more readable because of you. 	
\end{acknowledge}

%
%
%
%

\bibliographystyle{amsalpha}
\bibliography{bibliography}
\end{document}